\documentstyle[twoside,amstex]{article}
\input amssym.def
\input amssym.tex
\def\bc{\begin{center}}
\def\ec{\end{center}}

\begin{document}
\thispagestyle{empty} \vspace*{3 true cm} \pagestyle{myheadings}
\markboth {\hfill {\sl Nahid Ashrafi, Marjan Shebani and Huanyin
Chen}\hfill} {\hfill{\sl Rings Involving Idempotents, Units and
Nilpotent Elements}\hfill} \vspace*{-1.5 true cm} \bc{\large\bf
RINGS INVOLVING IDEMPOTENTS, UNITS AND\vskip2mm NILPOTENT
ELEMENTS}\ec

\vskip6mm
\bc{{\bf Nahid Ashrafi}\\[1mm]
Faculty of Mathematics, Statistics and Computer Science\\
Semnan University, Semnan, Iran\\
nashrafi@@semnan.ac.ir}\ec

\bc{{\bf Marjan
Shebani}\\[1mm]
Faculty of Mathematics, Statistics and Computer Science\\
Semnan University, Semnan, Iran\\
m.sheibani1@@gmail.com}\ec

\bc{{\bf Huanyin Chen$^*$}\\[1mm]
Department of Mathematics, Hangzhou Normal University\\
Hangzhou 310036, China\\
huanyinchen@@aliyun.com}\ec

\begin{abstract} We define the concepts of weakly precious ring and
precious ring which generalize the notions of a weakly clean ring
and a nil-clean ring. We obtain some fundamental properties of
these rings. We also consider certain subclasses of such rings,
and then offer new kinds of weakly clean rings and nil-clean
rings. Finally, we completely determine when a ring consists
entirely of very idempotents, units, and nilpotent elements.


\end{abstract}

\section{Introduction}

\vskip4mm Let $R$ be a ring with an identity. Idempotents, units
and nilpotent elements play important rules in ring theory. The
motivation of this paper is to investigate the structures of
various rings involving such special elements. An element $a$ in a
ring is called very idempotent if $a$ or $-a$ is an idempotent. An
element $a\in R$ is called (weakly) clean if there exists a (very)
idempotent $e\in R$ and a unit such that $a=e+u$~\cite{A}. An
element $a\in R$ is (weakly) nil-clean provided that there exists
a (very) idempotent $e\in R$ and a nilpotent $w\in R$ such that
$a=e+w$~\cite{B} and ~\cite{D}. These inspire us introduce two
concepts. We call an element $a\in R$ is (weakly) precious if
there exists a (very) idempotent $e\in R$, a unit $u\in R$ and a
nilpotent $w\in R$ such that $a=e+u+w$. A ring $R$ is called a
weakly clean (weakly precious, nil-clean, precious) ring if every
element in $R$ is weakly clean (weakly precious, nil-clean,
precious). Many fundamental properties about commutative weakly
clean rings were obtained in ~\cite{Ah} and ~\cite{A}, and that
weakly nil-clean rings were comprehensive studied by Breaz et al.
in ~\cite{B}.

In this paper, we shall explore the structures of these rings. In
Section 2, we prove that the direct product $R=\prod R_{i}$ of
rings $R_i$ is weakly precious if and only if each $R_{i}$ is
weakly precious and at most one is not precious. Furthermore, we
show that the precious property is invariant for any Morita
context. In Section 3, we are concern on weakly clean rings and
nil-clean rings. Let $R$ be a commutative ring with at most three
maximal ideals. If $2\in U(R)$ and $J(R)$ is nil, we prove that
$R$ is weakly clean. This provides a new type of weakly clean
rings. A ring $R$ is abelian if every idempotent is central. We
show that if $R$ is abelian then $M_n(R)$ is nil-clean if and only
if $R/J(R)$ is Boolean and $M_n(J(R))$ is nil. This extend the
main results of Breaz et al. ~\cite{BGDT} and that of Ko\c{s}an et
al.~\cite{KLZ}. In the last section, we investigate when a ring
consists entirely of very idempotent, units, and nilpotent
elements. We prove that a ring consists entirely of very
idempotents, units and nilpotent elements if and only if $R$ is
isomorphic to one of the following: a Boolean ring; ${\Bbb
Z}_3\bigoplus {\Bbb Z}_3$; ${\Bbb Z}_3\oplus B$ where $B$ is a
Boolean ring; local ring with a nil Jacobson radical;
$M_2\big({\Bbb Z}_2\big)$ or $M_2\big({\Bbb Z}_3\big)$; or the
ring of a Morita context with zero pairings where the underlying
rings are ${\Bbb Z}_2$ or ${\Bbb Z}_3$. The structure of such
rings is thereby completely determined.

Throughout, all rings are associative with an identity. $M_n(R)$
and $T_n(R)$ will denote the ring of all $n\times n$ full matrices
and triangular matrices over $R$, respectively. $J(R)$ and $P(R)$
stand for the Jacobson radical and prime radical of $R$. $Id(R)=\{
e\in R~|~e^2=e\in R\},-Id(R)=\{ e\in R~|~e^2=-e\in R\}$, $U(R)$ is
the set of all units in $R$, and $N(R)$ is the set of all
nilpotent elements in $R$.

\section{Weakly Previous Rings}

\vskip4mm We start this section by indicating that "weakly
cleanness" and "weakly preciousness", "nil-cleanness" and
"preciousness" are not the same for elements in a ring.

\vskip4mm \hspace{-1.8em} {\bf Example 2.1.}\ \ {\it $(1)$ Every
weakly clean element in a ring is weakly previous, but the
converse is not true.}\vspace{-.5mm}
\begin{enumerate}
\item [(2)] {\it Every nil-clean element in a ring is previous, but the converse is not true.}
\end{enumerate} \vspace{-.5mm} {\it Proof.}\ \ $(1)$ Obviously, every
weakly clean element in a ring is weakly previous. But the
converse is not true. Consider the matrix $A=\left(
\begin{array}{cc}
3&9\\
-7&-2
\end{array}
\right)\in M_2({\Bbb Z})$. By ~\cite[Theorem 4]{An}, $A$ is not
clean. Thus, $(A,-A)\in M_2({\Bbb Z})\times M_2({\Bbb Z})$ is not
weakly clean. But $(A,-A)$ is previous, as it has the previous
decomposition
$$\big(\left(
\begin{array}{cc}
1&0\\
7&0
\end{array}
\right),\left(
\begin{array}{cc}
1&0\\
6&0
\end{array}
\right)\big)+\big(\left(
\begin{array}{cc}
-1&0\\
-13&1
\end{array}
\right),\left(
\begin{array}{cc}
-1&0\\
0&-1
\end{array}
\right)\big)+\big(\left(
\begin{array}{cc}
3&9\\
-1&-3
\end{array}
\right),\left(
\begin{array}{cc}
-3&-9\\
1&3
\end{array}
\right)\big).$$ Thus, $(A,-A)$ is weakly precious.

$(2)$ Let $a\in R$ be nil-clean. Then there exists an idempotent
$e\in R$ and a nilpotent $w\in R$ such that $a=e+w$, and so
$a=(1-e)+(2e-1)+w$. As $(2e-1)^2=1$, we see that $a\in R$ is
precious. Thus, every nil-clean element in a ring is previous. The
converse is not true. For instance, $-1\in {\Bbb Z}_3$ is not
nil-clean, but it is precious.\hfill$\Box$

\vskip4mm Example 2.1 shows that $\{~\mbox{weakly clean
elements}~\}\subsetneq $ $\{~\mbox{weakly precious elements}\}$
and $\{~\mbox{nil-clean elements}~\}\subsetneq \{~\mbox{precious
elements}~\}.$ Though weakly precious rings are rich, but there
indeed exist rings which are not weakly precious. Since
$Id\big({\Bbb Z}\big)=\{ 0,1\}, U\big({\Bbb Z}\big)=\{ 1,-1\}$ and
${\Bbb Z}$ has no nonzero nilpotent element, we easily check that
$5\in {\Bbb Z}$ is not weakly precious. Therefore, the ring ${\Bbb
Z}$ of all integers is not weakly previous. The purpose of this
section is to investigate when a ring is weakly previous or
previous. Clearly, every homomorphic image of weakly precious
rings is weakly precious. Further, we derive

\vskip4mm \hspace{-1.8em} {\bf Lemma 2.2.}\ \ {\it Let $I$ be a
nil ideal of a ring $R$. Then $R$ is weakly precious if and only
if $R/I$ is weakly precious.} \vskip2mm\hspace{-1.8em} {\it
Proof.}\ \ Let $R$ be a weakly precious ring. Then $R/I$ is weakly
precious. Now assume that $R/I$ is weakly precious. Let $a\in R$.
Then $\overline{a}= a+I=\overline{e} +\overline{u}+\overline{w}$
or $\overline{a}= a+I=\overline{-e} +\overline{u}+\overline{w}$
for an idempotent $\overline{e}\in R/I$, a unit $\overline{u}\in
R/I$ and a nilpotent $\overline{w}\in R/I$. As $I$ is a nil ideal
of $R$, we easily check that $u$ is a unit element in $R$ and $w$
is a nilpotent element. Since every idempotent lifts modulo $I$,
we can find an idempotent $e\in R$ such that
$\overline{e}=\overline{f}$. Hence, $a= e+u+(w+c)$ or $a=
-e+u+(w+c)$ for some $c\in I$. Write $w^m=0~(m\geq 1)$. Then
$(w+c)^m\in I$, and so $w+c\in R$ is a nilpotent element.
Therefore $a\in R$ is weakly precious, as required.\hfill$\Box$

\vskip4mm \hspace{-1.8em} {\bf Theorem 2.3.}\ \ {\it Let $R$ be a
ring. Then $R$ is weakly precious if and only if $R[[x]]/(x^n)
(n\in {\Bbb N})$ is weakly precious.} \vskip2mm\hspace{-1.8em}
{\it Proof.}\ \ Clearly, $R[[x]]/(x^n)=\{ a_0+a_1x+\cdots
+a_{n-1}x^{n-1}~|~a_0,\cdots ,a_{n-1}\in R\}$. Let $\alpha :
R[[X]]/(x^n)\longrightarrow R$ be a morphism such that
$\alpha(f)=f(0)$. It is easy to check that $\alpha$ is an
$R$-epimorphism and $ker\alpha$ is a nil ideal of $R$, and
therefore the result follows from Lemma 2.2.\hfill$\Box$

\vskip4mm An element $a\in R$ is strongly nilpotent if for every
sequence $a_0,a_1,\cdots ,a_i,\cdots$ such that $a_0 =a$ and
$a_{i+1}\in a_iRa_i$, there exists an $n$ with $a_n=0$. The prime
radical (i.e., the intersection of all prime ideals) of a ring is
exactly the set of all its strongly nilpotent elements.

\vskip4mm \hspace{-1.8em} {\bf Lemma 2.4.}\ \ {\it Let $R$ be a
ring. Then the following are equivalent:}\vspace{-.5mm}
\begin{enumerate}
\item [(1)] {\it $R$ is weakly precious.}
\vspace{-.5mm} \item [(2)] {\it $R/P(R)$ is weakly
precious.}\end{enumerate} \vspace{-.5mm} {\it Proof.}\ \ This is
obvious by Lemma 2.2, as $P(R)$ is nil.\hfill$\Box$

\vskip4mm Recall that a ring $R$ is 2-primal provided that every
nilpotent element of $R$ is strongly nilpotent~\cite{MMZ}.

\vskip4mm \hspace{-1.8em} {\bf Theorem 2.5.}\ \ {\it Let $R$ be
2-primal. Then the following are equivalent:}\vspace{-.5mm}
\begin{enumerate}
\item [(1)] {\it $R$ is weakly precious;}
\vspace{-.5mm} \item [(2)] {\it $R$ is weakly clean.}
\vspace{-.5mm} \item [(3)] {\it $R/P(R)$ is weakly
clean.}\end{enumerate} \vspace{-.5mm} {\it Proof.}\ \
$(1)\Rightarrow (2)$ Let $a\in R$. As $R$ is weakly precious,
$a=e+u+w$ or $a=-e+u+w$ for an idempotent $e\in R$, a unit $u\in
R$ and a nilpotent $w\in R$. This shows that $a=
e+u\big(1+u^{-1}w)$ or $a=-e+u\big(1+u^{-1}w)$. As $R$ is a
2-primal ring, we get $w\in P(R)$. Since $P(R)$ is a nil ideal of
$R$, $1+u^{-1}w\in U(R)$, and therefore $R$ is weakly clean.

$(2)\Rightarrow (3)$ is clear.

$(3)\Rightarrow (1)$ Since $R/P(R)$ is weakly clean, it is weakly
precious. Therefore we complete the proof, by Lemma
2.4.\hfill$\Box$

\vskip4mm A ring $R$ is called nil-semicommutative  if $ab=0$ in
$R$ implies that $aRb=0$ for every $a, b \in N(R)$ (see
\cite{MMZ}). For instance, every semicommutative ring (i.e.,
$ab=0$ in $R$ implies that $aRb=0$) is nil-semicommutative.

\vskip4mm \hspace{-1.8em} {\bf Corollary 2.6.}\ \ {\it Let $R$ be
nil-semicommutative. Then the following are
equivalent:}\vspace{-.5mm}
\begin{enumerate}
\item [(1)] {\it $R$ is weakly precious;}
\vspace{-.5mm} \item [(2)] {\it $R$ is weakly clean.}
\vspace{-.5mm} \item [(3)] {\it $R/P(R)$ is weakly
clean.}\end{enumerate} \vspace{-.5mm} {\it Proof.}\ \  Let $R$ be
a nil-semicommutative ring. Then, by ~\cite[Lemma 2.7]{MMZ}, $R$
is 2-primal, so the result follows from Theorem 2.5.\hfill$\Box$

\vskip4mm A ring $R$ a right (left) quasi-duo ring if every
maximal right (left) ideal of $R$ is an ideal. For instance, local
rings, duo rings and weakly right (left) duo rings are all right
(left) quasi-duo rings. We now derive

\vskip4mm \hspace{-1.8em} {\bf Proposition 2.7.}\ \ {\it Let $R$
be a right (left) quasi-duo ring. Then the following are
equivalent:}\vspace{-.5mm}
\begin{enumerate}
\item [(1)] {\it $R$ is weakly precious.}
\vspace{-.5mm} \item [(2)] {\it $R$ is weakly clean.}
\vspace{-.5mm} \item [(3)] {\it $R/P(R)$ is weakly
clean.}\end{enumerate} \vspace{-.5mm} {\it Proof.}\ \
$(1)\Rightarrow (2)$ Let $a\in R$. Then there exists a very
idempotent $e\in R$, a unit $u\in R$ and a nilpotent element $w\in
R$ such that $a=e+u+w$. As $R$ is right (left) quasi-duo, it
follows from ~\cite[Lemma 2.3]{Yu} that $w\in J(R)$. Thus,
$a=e+u\big(1+u^{-1}w\big)$; hence, $a\in R$ is weakly clean.

$(2)\Rightarrow (3)$ This is obvious.

$(3)\Rightarrow (1)$ Clearly, $R/P(R)$ is weakly precious, and
therefore the result follows, by Lemma 2.4.\hfill$\Box$

\vskip4mm \hspace{-1.8em} {\bf Lemma 2.8.}\ \ {\it Let $R$ be
weakly precious and $S$ be precious. Then $R\oplus S$ is weakly
precious.} \vskip2mm\hspace{-1.8em} {\it Proof.}\ \ Set $A=R\oplus
S$. Let $(a,b)\in A$. Then there exists an idempotent $e\in R$, a
unit $u\in R$ and a nilpotent $v\in R$ such that $a=e+u+v$ or
$a=-e+u+v$.

Case I. $a=e+u+v$. Then we have an idempotent $f\in S$, a unit
$s\in S$ and a nilpotent $w\in S$ such that $b=f+v+w$. Thus,
$(a,b)=(e,f)+(u,s)+(v,w)$, where $(e,f)\in A$ is an idempotent,
$(u,s)\in A$ is a unit and $(v,w)\in A$ is nilpotent.

Case II. $a=-e+u+v$. Then we have an idempotent $f\in S$, a unit
$s\in S$ and a nilpotent $w\in S$ such that $-b=f+s+w$. Thus,
$(a,b)=-(e,f)+(u,-s)+(v,-w)$, where $(e,f)\in A$ is an idempotent,
$(u,-s)\in A$ is a unit and $(v,-w)\in A$ is nilpotent.

Therefore we conclude that $(a,b)$ is the sum of a very
idempotent, a unit and a nilpotent element in $A$, hence that
result.\hfill$\Box$

\vskip4mm \hspace{-1.8em} {\bf Theorem 2.9.}\ \ {\it Let
$\{R_{i}\}$ be a family of rings. Then the direct product $R=\prod
R_{i}$ of rings $R_i$ is weakly precious if and only if each
$R_{i}$ is weakly precious and at most one is not precious.}
\vskip2mm\hspace{-1.8em} {\it Proof.}\ \ $\Longrightarrow $
Obviously, each $R_i$ is weakly precious. Suppose $R_{i_1}$ and
$R_{i_2} (i_1\neq i_2)$ are not precious. Then there exist some
$x_{i_j}\in R_{i_j} (j=1,2)$ such that $x_{i_1}\in R_{i_1}$ and
$-x_{i_2}\in R_{i_2}$ are not precious. Choose $x=(x_i)$ where
$x_i=0$ whenever $i\neq i_j (j=1,2)$. Then $x\pm e$ is not the sum
of a unit and a nilpotent for all idempotents $e\in R$. This gives
a contradiction. Therefore each $R_{i}$ is a weakly precious and
at most one is not precious.

$\Longleftarrow$ Suppose that $R_{i_0}$ is weakly precious and all
the others $R_i$ are precious. Then $\prod\limits_{i\neq
i_{0}}R_{i_0}$ is precious. In light of Lemma 2.8, We conclude
that $R$ is weakly precious.\hfill$\Box$

\vskip4mm \hspace{-1.8em} {\bf Corollary 2.10.}\ \ {\it Let
$L=\prod\limits_{i\in I}R_{i}$ be the direct product of rings
$R_i\cong R$ and $|I|\geq 2$. Then $L$ is weakly precious if and
only if $R$ is precious if and only if $L$ is precious.}

\vskip4mm \hspace{-1.8em} {\bf Lemma 2.11.}\ \ {\it Let $e=e^2 \in
R$ be such that $eRe$ is weakly precious and $(1-e)R(1-e)$ is
precious. Then $R$ is weakly precious.} \vskip2mm\hspace{-1.8em}
{\it Proof.}  Let $e\in R$ be an idempotent, we have $R\cong
\left(
\begin{array}{cc}
eRe&eR(1-e)\\
(1-e)Re&(1-e)R(1-e)
\end{array}
\right)$. Now suppose that $A=\left(
\begin{array}{cc}
a&b\\
c&d
\end{array}
\right)$ be an element of $R$. As $eRe$ is weakly precious and
$(1-e)R(1-e)$ is precious, $a=f+u+w$ or $a=-f+u+w$ for some
idempotent $f\in eRe$, $u\in U(eRe)$ and a nilpotent $w\in\ eRe$.

Case I. Let $a=f+u+w$. Then $d-cu^{-1}b \in (1-e)R(1-e)$, and so
$d-cu^{-1}b=g+v+z$ for some idempotent $g$, unit $v$ and nilpotent
$z\in (1-e)R(1-e)$. Now $$A= \left(
\begin{array}{cc}
f+u+w&b\\
c&g+v+z+cu^{-1}b
\end{array}
\right) =\left(
\begin{array}{cc}
f&0\\
0&g
\end{array}
\right)+  \left(
\begin{array}{cc}
u&b\\
c&v+cu^{-1}b
\end{array}
\right)+ \left(
\begin{array}{cc}
w&0\\
0&z
\end{array}
\right).$$ It is clear that $\left(
\begin{array}{cc}
f&0\\
0&g
\end{array}
\right)$ is an idempotent element of $R$ and  $\left(
\begin{array}{cc}
w&0\\
0&z
\end{array}
\right) $ is a nilpotent element of $R$, so we need only to show
that $\left(
\begin{array}{cc}
u&b\\
c&v+cu^{-1}b
\end{array}
\right)$ is a unit of $R$. One easily checks that
$$\left(
\begin{array}{cc}
e&0\\
-cu^{-1}&1-e
\end{array}
\right)\left(
\begin{array}{cc}
u&b\\
c&v+cu^{-1}b
\end{array}
\right)\left(
\begin{array}{cc}
e&-u^{-1}\\
0&1-e
\end{array}
\right)\\
=\left(
\begin{array}{cc}
u&0\\
0&v
\end{array}
\right).$$ Hence, $\left(
\begin{array}{cc}
u&b\\
c&v+cu^{-1}b
\end{array}
\right)$ is invertible. Thus, $A$ is precious.

Case II. Let $a=-f+u+w$. Then $-a=f-u-w$. By the similar way of
Case I, we see that $-A$ is precious, as required.\hfill$\Box$

\vskip4mm A Morita context $(R,S,M,N,\psi,\varphi)$ consists of
two rings $R$ and $S$, two bimodules $_RN_S$ and $_SM_R$, and a
pair of bimodule homomorphisms $\psi : N\bigotimes\limits_{S}M\to
R$ and $\varphi: M\bigotimes\limits_{R} N\to S$ which satisfy the
following associativity: $\psi \big(n\bigotimes m\big)n'=n\varphi
\big(m\bigotimes n'\big)$ and $\varphi \big(m\bigotimes
n\big)m'=m\psi \big(n\bigotimes m'\big)$ for any $m,m'\in M,
n,n'\in N$. The ring $T=\{ \left(
\begin{array}{cc}
r&m\\
n&s \end{array} \right)~|~r\in R,s\in S,m\in M,n\in N\}$ is called
the ring of the Morita context $(R,S,M,N,\psi,\varphi)$.

\vskip4mm \hspace{-1.8em} {\bf Theorem 2.12.}\ \ {\it Let $T$ be
the ring of the Morita context $(R,S,M,N,\varphi,\phi)$. If $R$ is
weakly precious and $S$ is precious, then $T$ is weakly precious.}
\vskip2mm\hspace{-1.8em} {\it Proof.}\ \ Let $R$ be weakly
precious and $S$ be precious, and let $e=diag(1_{R},0)$. Since
$eTe\cong R$ and $(1_{T}-e)T(1_{T}-e)\cong S$, it follows by Lemma
2.11 that $T$ is a weakly precious ring, as asserted.\hfill$\Box$

\vskip4mm Many properties of weakly precious rings can be extended
to precious rings. For instance, $R$ is precious if and only if
$R[[x]]/(x^n) (n\in {\Bbb N})$ is precious. The direct product
$\prod R_{i}$ of rings $R_i$ is precious if and only if each
$R_{i}$ is precious. But the subdirect product of (weakly)
precious rings is not necessarily (weakly) precious. For instance,
${\Bbb Z}$ is a subdirect product of rings $\lbrace{\Bbb Z_n},
n\geq 2\rbrace$, where each ${\Bbb Z_n} (n\geq 2)$ is precious,
but ${\Bbb Z}$ is not.

\vskip4mm \hspace{-1.8em} {\bf Lemma 2.13.}\ \ {\it Let $e=e^2 \in
R$ be such that $eRe$ and $(1-e)R(1-e)$ are precious. Then $R$ is
precious.} \vskip2mm\hspace{-1.8em} {\it Proof.}\ \ Let $e\in R$
be an idempotent element, we have $R\cong \left(
\begin{array}{cc}
eRe&eR(1-e)\\
(1-e)Re&(1-e)R(1-e)
\end{array}
\right)$. Now suppose that  $A=\left(
\begin{array}{cc}
a&b\\
c&d
\end{array}
\right)$ be an element of $R$. Since $eRe$ and $(1-e)R(1-e)$ are
precious rings, $a=f+u+w$ for some idempotent $f\in eRe$, $u\in
U(eRe)$ and a nilpotent $w\in\ eRe$. Let $u^{-1}$ be the inverse
of $u$. Then $d-cu^{-1}b\in (1-e)R(1-e)$, so $d-cu^{-1}b=g+v+z$
for some idempotent $g$, unit $v$ and nilpotent $z\in
(1-e)R(1-e)$. Now $A= \left(
\begin{array}{cc}
f+u+w&b\\
c&g+v+z+cu^{-1}b
\end{array}
\right) =\left(
\begin{array}{cc}
f&0\\
0&g
\end{array}
\right)+\left(
\begin{array}{cc}
u&0\\
0&v+cu^{-1}b
\end{array}
\right)+ \left(
\begin{array}{cc}
w&0\\
0&z
\end{array}
\right)$. As in the proof of Lemma 2.11, we easily checks that
$\left(
\begin{array}{cc}
u&b\\
c&v+cu^{-1}b
\end{array}
\right)$ is a unit of $R$, and the the result follows.\hfill$\Box$

\vskip4mm \hspace{-1.8em} {\bf Theorem 2.14.}\ \ {\it Let $T$ be
the ring of the Morita context $(R,S,M,N,\varphi,\phi)$. If $R$
and $S$ are precious, then $T$ is precious.}
\vskip2mm\hspace{-1.8em} {\it Proof.}\ \ Let $R,S$ be precious
rings and let $e=diag(1_{R},0)$. Then $eTe\cong R$ and
$(1_{T}-e)T(1_{T}-e)\cong S$. By virtue of Lemma 2.13, $T$ is a
precious ring. \hfill$\Box$

\vskip4mm \hspace{-1.8em} {\bf Corollary 2.15.}\ \ {\it Let $R$ be
precious. Then $M_n(R)$ is precious.} \vskip2mm\hspace{-1.8em}
{\it Proof.}\ \ If $n=2$. Then the result follows by Theorem 2.14.
Suppose that the result holds for $n\leq k ( k\geq 2)$. Then $R$
and $M_k(R)$ are both precious. In light of Theorem 2.14, the ring
$\left(
\begin{array}{cc}
R&M\\
N&M_k(R)
\end{array}
\right)$ is precious, where $M=\{\left(
\begin{array}{ccc}
b_{1}&\cdots &b_k
\end{array}
\right)~|~b_1,\cdots ,b_k\in R\}$ and $N=\{\left(
\begin{array}{c}
c_{1}\\
\vdots\\
c_k
\end{array}
\right)~|~c_1,\cdots ,c_k\in R\}$. This completes the proof by
induction.\hfill$\Box$

\vskip4mm \hspace{-1.8em} {\bf Theorem 2.16.}\ \ {\it Let $R$ be a
ring. Then the following are equivalent:}\vspace{-.5mm}
\begin{enumerate}
\item [(1)] {\it $R$ is precious.}
\vspace{-.5mm}
\item [(2)] {\it $T_n(R)$ is precious for all $n\in {\Bbb N}$.}
\vspace{-.5mm}
\item [(3)] {\it
$T_n(R)$ is precious for some $n\in {\Bbb N}$.} \vspace{-.5mm}
\item [(4)] {\it
$T_n(R)$ is weakly precious for some $n\geq 2$.}
\end{enumerate} \vspace{-.5mm} {\it Proof.}\ \ $(1)\Rightarrow
(2)$  The result holds for $n=2$ by Theorem 2.14. Assume that the
result holds for $n\leq k$ $(k\geq 2)$. Let $n=k+1$. Then
$T_n(R)\cong \left(
\begin{array}{cc} R&M\\
0&T_{k}(R) \end{array} \right)$, where $M =\{ ( c_{1}, \cdots ,
c_{k})\mid c_1,\cdots ,c_{k}\in R\}$. In light of Theorem 2.14,
$T_{k+1}(R)$ is precious. Then, proving $(2)$, by induction.

$(2)\Rightarrow (3)$ is trivial.

$(3)\Rightarrow (1)$ This is obvious, as a triangular matrix over
$R$ is an idempotent (unit, nilpotent matrix) if and only if every
its diagonal entry is an idempotent (unit, nilpotent matrix).

$(2)\Rightarrow (4)$ is trivial.

$(4)\Rightarrow (1)$ Let $a\in R$. Choose $A=\left(
\begin{array}{cccccc}
a&&&&&\\
&-a&&&&\\
&&0&&&\\
&&&\ddots&&\\
&&&&&0
\end{array}
\right)\in T_n(R)$. By hypothesis, we can find an idempotent
$\left(
\begin{array}{ccccc}
e_1&&&&\\
&e_2&&&\\
&&\ddots&&\\
&&&&e_n
\end{array}
\right)$, a unit $\left(
\begin{array}{ccccc}
u_1&&&&\\
&u_2&&&\\
&&\ddots&&\\
&&&&u_n
\end{array}
\right)$ and a nilpotent $\left(
\begin{array}{ccccc}
w_1&&&&\\
&w_2&&&\\
&&\ddots&\\
&&&&w_n
\end{array}
\right)$ such that
$$A=\pm
\left(
\begin{array}{ccccc}
e_1&&&&\\
&e_2&&&\\
&&\ddots&&\\
&&&&e_n
\end{array}
\right)+\left(
\begin{array}{ccccc}
u_1&&&&\\
&u_2&&&\\
&&\ddots&&\\
&&&&u_n
\end{array}
\right)+\left(
\begin{array}{ccccc}
w_1&&&&\\
&w_2&&&\\
&&\ddots&&\\
&&&&w_n
\end{array}
\right).$$ It follows that $a=e_1+u_1+w_1$ or $a=e_2-u_2-w_2$.
Clearly, $e_1,e_2$ are idempotents, $u_1, u_2$ are units and
$w_1,w_2$ are nilpotent. Therefore proving $(1)$.\hfill$\Box$

\vskip4mm A ring $R$ is weakly periodic provided that for any
$a\in R$ there exists some $p=p^m(m\geq 2)$ such that $a-p\in R$
is nilpotent. For instance, every periodic ring is weakly
periodic.

\vskip4mm \hspace{-1.8em} {\bf Corollary 2.17.}\ \ {\it Let $R$ be
a weakly periodic ring. Then $M_n(R)$ and $T_n(R)$ are precious
for all $n\in {\Bbb N}$.} \vskip2mm\hspace{-1.8em} {\it Proof.}\ \
For any $a\in R$, there exists a $p=p^{k+1} (k\in {\Bbb N})$ and a
nilpotent $w\in R$ such that $a=p+w$. If $k=1$, then $p\in R$ is
an idempotent, and so $a\in R$ is precious. Suppose that $k\geq
2$. Set $e=1-p^k, u=p-1+p^k$. Then $e=e^2\in
R,u^{-1}=p^{k-1}-1+p^k$, and that $p=e+u$. Thus, $a=e+u+w$, and
then $R$ is precious. Therefore we complete the proof, by
Corollary 2.15 and Theorem 2.16.\hfill$\Box$

\vskip4mm Let $R$ be a ring and $M$ an $R$-$R$-bimodule. The
trivial extension of $R$ by $M$,
$$R\propto M=\{ \left(
\begin{array}{cc}
r&m\\
&r
\end{array}
\right)~|~r\in R,m\in M\}$$ is (weakly) precious if and only if
$R$ is (weakly) precious.

\section{Certain Subclasses}

\vskip4mm Weakly clean rings and nil-clean rings forms main types
of subclasses of weakly precious rings and precious rings,
respectively. The purpose of this section is to off new types of
such rings. In ~\cite{A}, Anderson and Camillo proved that if a
ring $R$ has at most two maximal ideals and $2\in U(R)$ then $R$
is weakly clean. We extend this result as follows.

\vskip4mm \hspace{-1.8em} {\bf Theorem 3.1.}\ \ {\it Let $R$ be a
commutative ring with at most three maximal ideals. If $2\in U(R)$
and $J(R)$ is nil, then $R$ is weakly clean.}
\vskip2mm\hspace{-1.8em} {\it Proof.}\ \ Case I. $R$ has only one
maximal ideal. Then $R$ is local, and so it is clean.

Case II. $R$ has only two maximal ideals. Then $R$ is weakly
clean, by ~\cite[Proposition 16]{A}.

Case III. $R$ has three maximal ideals $M_1,M_2$ and $M_2$. Let
$a\in R$. If $a\not\in M_1,M_2,M_3$, then $aR=R$; hence, $a\in
U(R)$. So we may assume that $a\in M_1$. If $1-a\not\in
M_1,M_2,M_3$, then $1-a\in U(R)$, and so $a\in R$ is clean. If
$1-a\in M_1$, then $1=a+(1-a)\in M_1$, a contradiction. Thus,
$1-a\in M_2$ or $1-a\in M_3$. If $1+a\not\in M_1,M_2,M_3$, then
$1+a\in U(R)$; hence, $a\in R$ is weakly clean. If $1+a\in M_1$,
then $1=(1+a)-a\in M_1$, a contradiction. Thus, $1+a\in M_2$ or
$1+a\in M_3$. There are only two possible cases: $1-a\in M_2,
1+a\in M_3$ or $1-a\in M_3, 1+a\in M_2$, as $2\in U(R)$. Thus, we
may assume that $1-a\in M_2, 1+a\in M_3$. Hence, $a(1-a)(1+a)\in
M_1M_2M_3\subseteq J(R)$. Thus, $\overline{a}=\overline{a}^3\in
R/J(R)$. Set $\overline{e}=\overline{1-a^2}$ and
$\overline{u}=\overline{a^2+a-1}$. Then $\overline{e}\in R/J(R)$
is an idempotent and $\overline{u}^2=\overline{1}$. Further, we
have $\overline{a}=\overline{e}+\overline{u}$. By hypothesis,
$J(R)$ is nil, and then every unit and every idempotent lift
modulo $J(R)$. So we may assume that $e\in R$ is an idempotent and
$u\in U(R)$. Set $w:=a-e-u$. Then $a=e+u+w$ where $w\in J(R)$.
Clearly, $u+w=u(1+u^{-1}w)\in U(R)$, and so $a\in R$ is clean.
Therefore $R$ is weakly clean.\hfill$\Box$

\vskip4mm \hspace{-1.8em} {\bf Example 3.2.}\ \ {\it Let
$R=k[[x,y,z]]$ where $k$ is a field with $char(k)\neq 2$. Let
$S=R-(x)\cup (y)\cup (z)$. Then the ring $R_S/J^2(R_S)$ is weakly
clean.} \vskip2mm\hspace{-1.8em} {\it Proof.}\ \ Since
$k[[x,y,z]]/(x)\cong k[[y,z]]$ is an integral domain, we see that
$(x)$ is a prime ideal of $k[[x,y,z]]$. Likewise, $(y)$ and $(z)$
are prime ideals of $R$. Let $S=R-(x)\bigcup (y)\bigcup (z)$. Then
$S$ is a multiplicative closed subset of $R$. Let $P$ be a maximal
ideal of $R_S$. Then we have an ideal $Q$ of $R$ such that $P=Q_S$
such that $Q\bigcap S=\emptyset$. Thus, $Q\subseteq (x)\bigcup
(y)\bigcup (z)$. Assume that $Q\nsubseteq (x), Q\nsubseteq (y)$
and $Q\nsubseteq (z)$. Then we have some $b,c,d\in Q$, but
$b\not\in (x)$, $c\not\in (y)$ and $d\not\in (z)$. Choose
$a=b+c+d$. Then $a\in Q$. If $a\in (x)$, then $c+d\in (x)$. If
$c\not\in (x)$, then $c\in (z)$. This implies that $d\in (x)$.
Hence, $c\in (z)\bigcap (x)=0$. This gives a contradiction. If
$c\in (x)$, then $d\in (x)$; hence that $b=a-(c+d)\in (x)$, a
contradiction. Hence, $a\not\in (x)$. Likewise, $a\not\in (y)$ and
$a\not\in (z)$. Thus, $a\not\in (x)\bigcup (y)\bigcup (z)$, a
contradiction. We infer that $Q\subseteq (x)$, or $Q\subseteq (y)$
or $Q\subseteq (z)$. Hence, $Q_S\subseteq (x)_S$, or $Q_S\subseteq
(y)_S$, or $Q_S\subseteq (z)_S$. By the maximality of $P$, we get
$P=(x)_S$, or $(y)_S$, or $(z)_S$. Thus, $R_S$ has exactly three
maximal ideals $(x)_S$, $(y)_S$ and $(z)_S$. Therefore $R$ has at
most three maximal ideals. Since $char(k)\neq 2$, we see that
$2\in U(R_S)$.

Set $A=R_S/J^2(R_S)$. Then $A$ has at most three maximal ideals
and $2\in U(A)$. If $\overline{x}\in J(A)$, then
$\overline{1-xr}\in U(A)$ for any $r\in R_S$. Hence, $1-xr\in
U(R_S)$. This implies that $x\in J(R_S)$, and so
$\overline{x^2}=\overline{0}$. That is, $\overline{x}$ is
nilpotent. So, $J(A)$ is nil. Therefore we complete the proof, in
terms of Theorem 3.1.\hfill$\Box$

\vskip4mm \hspace{-1.8em} {\bf Example 3.3.}\ \ {\it Let $R=\{
\frac{m}{n}~|~(m,n)=1, m,n\in {\Bbb Z}~\mbox{and}~3,5,7~\nmid
n\}$. Then the ring $R/J^2(R)$ is weakly clean.}
\vskip2mm\hspace{-1.8em} {\it Proof.}\ \ Let $M$ be an ideal of
$R$ such that $3R\subsetneq M\subseteq R$. Choose $\frac{m}{n}\in
M$ while $\frac{m}{n}\not\in 3R$. Then $3\nmid m$, and so
$(3,m)=1$. So $3k+lm=1$ for some $k,l\in {\Bbb Z}$. This shows
that $\frac{1}{1}=3\cdot \frac{k}{1}+l\frac{m}{n}\cdot
\frac{n}{1}\in M$, i.e., $M=R$. Thus, $3R$ is a maximal ideal of
$R$. Likewise, $5R$ and $7R$ are maximal ideals of $R$. For any
$\frac{m}{n}\in 3R\bigcap 5R\bigcap 7R$ and $\frac{a}{b}\in R$,
then $\frac{1}{1}-\frac{m}{n}\frac{a}{b}=\frac{nb-ma}{nb}$. Write
$\frac{m}{n}=\frac{3s}{t}$. Then $3sn=mt$, and so $3~|~mt$. Since
$3\nmid t$, we get $3~|~m$. Obviously, $3\nmid nb$; hence, $3\nmid
(nb-ma)$. Similarly, $5,7\nmid (nb-ma)$. It follows that
$\frac{nb}{nb-ma}\in U(R)$. We infer that $\frac{m}{n}\in J(R)$.
Therefore $3R\bigcap 5R\bigcap 7R\subseteq J(R)$. Let $M$ be a
maximal ideal of $R$ and $M\neq 3R, 5R, 7R$. Then $3R+M=R, 5R+M=R$
and $7R+M=R$. Thus, $R=(3R+M)(5R+M)(7R+M)\subseteq 3R\bigcap
5R\bigcap 7R+M=J(R)+M\subseteq M$, hence, $R=M$, an absurd. We
infer that $R$ is a commutative ring with exactly three maximal
ideals. Obviously $2\in R$ is invertible. Therefore $A:=R/J^2(R)$
is a commutative ring with exactly three maximal ideals. Obviously
$2\in A$ is invertible. As in the proof of Example 3.2, $A$ has
the nil Jacobson radical. We conclude that $A$ is weakly clean, by
Theorem 3.1.\hfill$\Box$

\vskip4mm In ~\cite[Question 3]{D}, Diesl asked: Let $R$ be a nil
clean ring, and let $n$ be a positive integer. Is $M_n(R)$ nil
clean?  In ~\cite[Theorem 3]{BGDT}, Breaz et al. proved that their
main theorem: for a field $K$, $M_n(K)$ is nil-clean if and only
if $K\cong {\Bbb Z}_2$. They also asked if this result could be
extended to division rings. As a main result in ~\cite{KLZ},
Ko\c{s}an et al. gave a positive answer to this problem. They
showed that the preceding equivalence holds for any division ring.
We shall extend ~\cite[Theorem 3]{BGDT} and ~\cite[Theorem 3]{KLZ}
to an arbitrary abelian ring.

\vskip4mm Recall that a ring $R$ is an exchange ring if for every
$a\in R$ there exists an idempotent $e\in aR$ such that $1-e\in
(1-a)R$. Clearly, every nil-clean ring is an exchange ring.

\vskip4mm \hspace{-1.8em} {\bf Lemma 3.4.}\ \ {\it Let $R$ be an
abelian exchange ring, and let $x\in R$. Then $RxR=R$ if and only
if $x\in U(R)$.} \vskip2mm\hspace{-1.8em} {\it Proof.}\ \ If $x\in
U(R)$, then $RxR=R$. Conversely, assume that $RxR=R$. As in the
proof of ~\cite[Proposition 17.1.9]{CH}, there exists an
idempotent $e\in R$ such that $e\in xR$ such that $ReR=R$. This
implies that $e=1$. Write $xy=1$. Then $yx=y(xy)x=(yx)^2$. Hence,
$yx=y(yx)x$. Therefore $1=x(yx)y=xy(yx)xy=yx$, and so $x\in U(R)$.
This completes the proof.\hfill$\Box$

\vskip4mm Set $J^*(R)=\bigcap \{ P~|~P~\mbox{is a maximal ideal
of}~R\}.$ We will see that $J(R)\subseteq J^*(R)$. In general,
they are not the same. For instance, $J(R)=0$ and $J^*(R)=\{ x\in
R~|~dim_F(xV)<\infty\}$, where $R=End_F(V)$ and $V$ is an
infinite-dimensional vector space over a field $F$.

\vskip4mm \hspace{-1.8em} {\bf Lemma 3.5.}\ \ {\it Let $R$ be an
abelian exchange ring. Then $J^*(R)=J(R)$.}
\vskip2mm\hspace{-1.8em} {\it Proof.}\ \ Let $M$ be a maximal
ideal of $R$. If $J(R)\nsubseteq M$, then $J(R)+M=R$. Write
$x+y=1$ with $x\in J(R),y\in M$. Then $y=1-x\in U(R)$, an absurd.
Hence, $J(R)\subseteq M$. This implies that $J(R)\subseteq
J^*(R)$. Let $x\in J^*(R)$, and let $r\in R$. If $R(1-xr)R\neq R$,
then we can find a maximal ideal $M$ of $R$ such that
$R(1-xr)R\subseteq M$, and so $1-xr\in M$. It follows that
$1=xr+(1-xr)\in M$, which is imposable. Therefore $R(1-xr)R=R$. In
light of Lemma 3.4, $1-xr\in U(R)$, and then $x\in J(R)$. This
completes the proof.\hfill$\Box$

\vskip4mm \hspace{-1.8em} {\bf Lemma 3.6.}\ \ {\it Let $R$ be a
ring with no non-trivial idempotents, and let $n\in {\Bbb N}$.
Then the following are equivalent:}\vspace{-.5mm}
\begin{enumerate}
\item [(1)]{\it $M_n(R)$ is nil-clean.}
\item [(2)]{\it $R/J(R)\cong {\Bbb Z}_2$ and $M_n(J(R))$ is nil.}
\end{enumerate} \vspace{-.5mm} {\it Proof.}\ \ $(1)\Rightarrow (2)$ In view of ~\cite[Proposition 3.16]{D}, $J(M_n(R))$ is
nil, and then so is $M_n(J(R))$.

Let $a\in R$. By hypothesis, $M_n(R)$ is nil-clean. If $n=1$, then
$R$ is nil-clean. Hence, $a\in N(R)$ or $a-1\in N(R)$. This shows
that $a\in U(R)$ or $1-a\in U(R)$, and so $R$ is local. That is,
$R/J(R)$ is a division ring. As $R/J(R)$ is nil-clean, it follows
from ~\cite[Theorem 3]{BGDT} that $R/J(R)\cong {\Bbb Z}_2$. We now
assume that $n\geq 2$. Then there exists an idempotent $E\in
M_n(R)$ and a nilpotent $W\in GL_n(R)$ such that $I_n+\left(
\begin{array}{cccc}
a&&\\
&0&\\
&&\ddots&\\
&&&0
\end{array}
\right)=E+W$. Set $U=-I_n+W$. Then $U\in GL_n(R)$. Hence,
$$U^{-1}\left(
\begin{array}{cccc}
a&&\\
&0&\\
&&\ddots&\\
&&&0
\end{array}
\right)=U^{-1}E+I_n=\big(U^{-1}EU\big)U^{-1}+I_n.$$ Set
$F=U^{-1}EU$. Then $F=F^2\in M_n(R)$, and that
$$(I_n-F)U^{-1}\left(
\begin{array}{cccc}
a&&\\
&0&\\
&&\ddots&\\
&&&0
\end{array}
\right)=I_n-F.$$ Write $I_n-F=\left(
\begin{array}{cccc}
e&0&\\
*&0&\\
\vdots&&\ddots&\\
*&0&&0
\end{array}
\right).$ By hypothesis, $e=0$ or $1$. If $e=0$, then $I_n-F=0$,
and so $E=I_n$. This shows that $\left(
\begin{array}{cccc}
a&&\\
&0&\\
&&\ddots&\\
&&&0
\end{array}
\right)=W$ is nilpotent; hence that $a\in R$ is nilpotent. Thus,
$1-a\in U(R)$.

If $e=1$, then $F=\left(
\begin{array}{cccc}
0&0&\\
*&1&\\
\vdots&&\ddots&\\
*&0&&1
\end{array}
\right).$ Write $U^{-1}=\left(
\begin{array}{cc}
\alpha&\beta\\
\gamma&\delta
\end{array}
\right),$ where $\alpha\in R, \beta\in M_{1\times (n-1)}(R),$
$\gamma\in M_{(n-1)\times 1}(R)$ and $\delta\in M_{(n-1)\times
(n-1)}(R)$. Then $$\left(
\begin{array}{cc}
\alpha&\beta\\
\gamma&\delta
\end{array}
\right)\left(
\begin{array}{cccc}
a&&\\
&0&\\
&&\ddots&\\
&&&0
\end{array}
\right)=\left(
\begin{array}{cc}
0&0\\
x&I_{n-1}
\end{array}
\right)\left(
\begin{array}{cc}
\alpha&\beta\\
\gamma&\delta
\end{array}
\right)+I_n,$$ where $x\in M_{(n-1)\times 1}(R)$. Thus, we get
$$\begin{array}{c}
\alpha a=1, \gamma a=x\alpha+\gamma, 0=x\beta+\delta+I_{n-1}.
\end{array}$$ One easily checks that
$$\left(
\begin{array}{cc}
1&\beta\\
0&I_{n-1}
\end{array}
\right)\left(
\begin{array}{cc}
1&0\\
x&I_{n-1}
\end{array}
\right)U^{-1}\left(
\begin{array}{cc}
1&0\\
\gamma a&I_{n-1}
\end{array}
\right)=\left(
\begin{array}{cc}
\alpha+\beta\gamma a&0\\
0&-I_{n-1}
\end{array}
\right).$$ This implies that $u:=\alpha+\beta\gamma a\in U(R)$.
Hence, $\alpha=u-\beta\gamma a$. It follows from $\alpha a=1$ that
$(u-\beta\gamma a)a=1$. Since $R$ has only trivial idempotents, we
get $a(u-\beta\gamma a)=1$, and so $a\in U(R)$. This shows that
$a\in U(R)$ or $1-a\in U(R)$. Therefore $R$ is local, and then
$R/J(R)$ is a division ring. Since $M_n(R)$ is nil-clean, we see
that so is $M_n(R/J(R))$. In light of ~\cite[Theorem 3]{BGDT},
$R/J(R)\cong {\Bbb Z}_2$, as desired.

$(2)\Rightarrow (1)$ By virtue of ~\cite[Theorem 3]{BGDT},
$M_n(R/J(R))$ is
 nil-clean. Since $M_n(R)/J(M_n(R))\cong M_n(R/J(R))$ and
$J\big(M_n(R)\big)=M_n(J(R))$ is nil, it follows from ~\cite[Lemma
4]{BGDT} that $M_n(R)$ is nil-clean, as asserted. \hfill$\Box$

\vskip4mm \hspace{-1.8em} {\bf Example 3.7.}\ \ Let $K$ be a
field, and let $R=K[x,y]/(x,y)^2$. Then $M_n(R)$ is nil-clean if
and only if $K\cong {\Bbb Z}_2$. Clearly, $J(R)=(x,y)/(x,y)^2$,
and so $R/J(R)\cong K$. Thus, $R$ is a local ring with a nilpotent
Jacobson radical. Hence, $R$ has no non-trivial idempotents. Thus,
we are done by Lemma 3.6.

\vskip4mm We are now ready to prove:

\vskip4mm \hspace{-1.8em} {\bf Theorem 3.8.}\ \ {\it Let $R$ be
abelian, and let $n\in {\Bbb N}$. Then the following are
equivalent:}\vspace{-.5mm}
\begin{enumerate}
\item [(1)]{\it $M_n(R)$ is nil-clean.}
\item [(2)]{\it $R/J(R)$ is Boolean and $M_n(J(R))$ is nil.}
\end{enumerate} \vspace{-.5mm} {\it Proof.}\ \ $(1)\Rightarrow (2)$ Clearly, $M_n(J(R))$ is nil. Let $M$ be a maximal ideal of
$R$, and let $\varphi_M: R\to R/M.$ Since $M_n(R)$ is nil-clean,
then so is $M_n(R/M)$. Hence, $R/M$ is an exchange ring with all
idempotents central. In view of ~\cite[Lemma 17.2.5]{CH}, $R/M$ is
local, and so $R/M$ has only trivial idempotents. It follows from
Lemma 3.6 that $R/M/J(R/M)\cong {\Bbb Z}_2$. Write $J(R/M)=K/M$.
Then $K$ is a maximal ideal of $R$, and that $M\subseteq K$. This
implies that $M=K$; hence, $R/M\cong {\Bbb Z}_2$. Construct a map
$\varphi_M: R/J^*(R)\to R/M, r+J^*(R)\mapsto r+M$. Here, $J^*(R)$
is the intersection of all maximal two-sided ideal of $R$. Then
$\bigcap\limits_{M}Ker\varphi_M=\bigcap\limits_{M}\{
r+J^*(R)~|~r\in M\}=0$. Therefore $R/J^*(R)$ is isomorphic to a
subdirect product of some ${\Bbb Z}_2$. Hence, $R/J^*(R)$ is
Boolean. In light of Lemma 3.5, $R/J(R)$ is Boolean, as desired.

$(2)\Rightarrow (1)$ Since $R/J(R)$ is Boolean, it follows by
~\cite[Corollary 6]{BGDT} that $M_n(R/J(R))$ is nil-clean. That
is, $M_n(R)/J(M_n(R))$ is nil-clean. But $J(M_n(R))=M_n(J(R))$ is
nil. Therefore we complete the proof, in terms of Lemma
3.6.\hfill$\Box$

\vskip4mm We note that the "$(2)\Rightarrow (1)$" in Theorem 3.8
always holds, but "abelian" condition is necessary in
"$(1)\Rightarrow (2)$". Let $R=M_n({\Bbb Z}_2) (n\geq 2)$. Then
$R$ is nil-clean. But $R/J(R)$ is not Boolean. Here, $R$ is not
abelian.

\vskip4mm \hspace{-1.8em} {\bf Corollary 3.9.}\ \ {\it Let $R$ be
commutative, and let $n\in {\Bbb N}$. Then the following are
equivalent:}\vspace{-.5mm}
\begin{enumerate}
\item [(1)]{\it $M_n(R)$ is nil-clean.}
\item [(2)]{\it $R/J(R)$ is Boolean and $J(R)$ is nil.}
\item [(3)]{\it For any $a\in R$, $a-a^2\in R$ is nilpotent.}
\end{enumerate} \vspace{-.5mm} {\it Proof.}\ \ $(1)\Rightarrow (3)$ Let $a\in R$. In view of Theorem 3.8,
$a-a^2\in J(R)$. Since $R$ is commutative, we see that $J(R)$ is
nil if and only if $J(M_n(R))$ is nil. Therefore $a-a^2\in R$ is
nilpotent.

$(3)\Rightarrow (2)$ Clearly, $R/J(R)$ is Boolean. For any $a\in
J(R)$, we have $(a-a^2)^n=0$ for some $n\geq 1$. Hence,
$a^n(1-a)^n=0$, and so $a^n=0$. This implies that $J(R)$ is nil.

$(2)\Rightarrow (1)$ As $R$ is commutative, we see that
$M_n(J(R))$ is nil. This completes the proof, by Theorem 3.8.
\hfill$\Box$

\vskip4mm Furthermore, we observe that the converse of
~\cite[Corollary 7]{BGDT} is true as the following shows.

\vskip4mm \hspace{-1.8em} {\bf Corollary 3.10.}\ \ {\it A
commutative ring $R$ is nil-clean if and only if $M_n(R)$ is
nil-clean.} \vskip2mm\hspace{-1.8em} {\it Proof.}\ \ One direction
is obvious by ~\cite[Corollary 7]{BGDT}. Suppose that $M_n(R)$ is
nil-clean. In view of Corollary 3.9, $R/J(R)\cong {\Bbb Z}_2$ is
nil-clean, and that $J(R)$ is nil. Therefore $R$ is nil-clean, by
~\cite[Lemma 4]{BGDT}.\hfill$\Box$

\vskip4mm \hspace{-1.8em} {\bf Corollary 3.11.}\ \ Let $m,n\in
{\Bbb N}$. Then $M_n\big({\Bbb Z}_m\big)$ is nil-clean if and only
if $m=2^r$ for some $r\in {\Bbb N}$. Write $m=p_1^{r_1}\cdots
p_s^{r_s} (p_1,\cdots ,p_s~\mbox{are distinct primes}, r_1,\cdots
,r_s\in {\Bbb N}$). Then $Z_m\cong Z_{p_1^{r_1}}\oplus \cdots
\oplus Z_{p_m^{r_s}}$. In light of Corollary 3.10, $M_n\big({\Bbb
Z}_m\big)$ is nil-clean if and only if ${\Bbb Z}_m$ is nil-clean,
if and only if $s=1$ and $Z_{p_1^{r_1}}$ is nil-clean. Therefore
we are done by Lemma 3.6.

\section{A Special Case}

\vskip4mm A natural problem is asked when a ring consists entirely
of very idempotents, units, and nilpotent elements. We will extend
the study of the rings consisting entirely of some special
elements in ~\cite{I}, and explore such type rings. Surprisingly,
our case will be involved in both Boolean rings and the ring
${\Bbb Z}_3$ of integers modulo $3$. Their structures will be
thereby completely determined. The following is a generalization
of ~\cite[Corollary 2.29]{Ah} which is for a commutative ring.

\vskip4mm \hspace{-1.8em} {\bf Lemma 4.1.}\ \ {\it Let $R$ be a
ring. Then $R=U(R)\bigcup Id(R)\bigcup -Id(R)$ if and only if $R$
is isomorphic to one of the following:}\vspace{-.5mm}
\begin{enumerate}
\item [(1)] {\it a Boolean ring;} \vspace{-.5mm}
\item [(2)] {\it a division ring;}
\vspace{-.5mm}
\item [(3)] {\it ${\Bbb Z}_3\oplus {\Bbb Z}_3$; or} \vspace{-.5mm}
\item [(4)] {\it
${\Bbb Z}_3\oplus B$ where $B$ is a Boolean ring.}
\end{enumerate} \vspace{-.5mm} {\it Proof.}\ \ $\Longrightarrow $
It is easy to check that $R$ is reduced; hence, it is abelian.

Case I. $R$ is indecomposable. Then $R$ is a division ring.

Case II. $R$ is decomposable. Then $R=A\oplus B$ where $A,B\neq
0$. If $0\neq x\in A$, then $(x,0)\in R$ is a very idempotent.
Hence, $x\in R$ is very idempotent. Hence, $A=Id(A)\bigcup
-Id(A)$. Likewise, $B=Id(B)\bigcup -Id(B)$. In view of
~\cite[Theorem 1.12]{Ah}, $A$ and $B$ are isomorphic to one of the
following:
\begin{enumerate}
\item [(1)] {\it ${\Bbb Z}_3$;}
\vspace{-.5mm} \item [(2)] {\it a Boolean ring;} \vspace{-.5mm}
\item [(3)] {\it ${\Bbb Z}_3\oplus B$ where $B$ is a Boolean ring.}
\end{enumerate}

Thus, $R$ is isomorphic to one of the following:
\begin{enumerate}
\item [(a)] {\it ${\Bbb Z}_3\oplus {\Bbb Z}_3$;}
\vspace{-.5mm} \item [(b)] {\it ${\Bbb Z}_3\oplus B$ where $B$ is
a Boolean ring;} \vspace{-.5mm}
\item [(c)] {\it ${\Bbb Z}_3\oplus {\Bbb Z}_3\oplus B$ where $B$ is
a Boolean ring;} \vspace{-.5mm}
\item [(d)] {\it a Boolean ring;}\end{enumerate}

Case (c). $(1,-1,0)\not\in U(R)\bigcup Id(R)\bigcup -Id(R)$, an
absurd.

Therefore we conclude that $R$ is one of Cases (a), (b) and (d),
as desired.

$\Longleftarrow $ Case (1). $R=Id(R)$. Case (2). $R=U(R)\bigcup
Id(R)$. Case (3). $U(R)=\{ (1,1),(1,-1),$ $(-1,1),$ $(-1,-1)$,
$Id(R)=\{ (0,0),(0,1),(1,0)\}$ and $-Id(R)=\{ (0,0),(0,-1),$ $
(-1,0),(-1,-1)\}$. Thus, $R=U(R)\bigcup Id(R)\bigcup -Id(R)$. Case
(4). $Id(R)=\{ (0,x),(1,x) ~|~x\in B\}$ and $-Id(R)=\{
(0,x),(-1,x)~|~x\in B$. Therefore $R=Id(R)\bigcup -Id(R)$, as
desired.\hfill$\Box$

\vskip4mm \hspace{-1.8em} {\bf Lemma 4.2.}\ \ {\it Let $R$ be a
decomposable ring. Then $R$ consists entirely of very idempotents,
units, and nilpotent elements if and only if $R$ is isomorphic to
one of the following:}\vspace{-.5mm}
\begin{enumerate}
\item [(1)] {\it a Boolean ring;}
\vspace{-.5mm}
\item [(2)] {\it ${\Bbb Z}_3\times {\Bbb Z}_3$;}
\vspace{-.5mm}
\item [(3)] {\it
${\Bbb Z}_3\oplus B$ where $B$ is a Boolean ring.}
\end{enumerate} \vspace{-.5mm} {\it Proof.}\ \ $\Longrightarrow$ Write $R=A\oplus B$
with $A,B\neq 0$. Then $A$ and $B$ are rings that consisting
entirely of very idempotents, units, and nilpotent elements. If
$0\neq x\in N(A)$, then $(x,1)\not\in Id(R)\bigcup -Id(R)\bigcup
U(R)\bigcup N(R)$. This shows that $A=U(A)\bigcup Id(A)\bigcup
-Id(A)$. Likewise, $B=U(B)\bigcup Id(B)\bigcup -Id(B)$. In light
of Lemma 4.1, $R$ is one of the following:
\begin{enumerate}
\item [(a)] {\it a Boolean ring;} \vspace{-.5mm}
\item [(b)] {\it $B\oplus D$ where $B$ is a Boolean ring and $D$ is a division ring;}
\vspace{-.5mm}
\item [(c)] {\it ${\Bbb Z}_3\oplus {\Bbb Z}_3\oplus B$ where $B$ is
a Boolean ring;} \vspace{-.5mm}
\item [(d)] {\it ${\Bbb Z}_3\oplus B$ where $B$ is
a Boolean ring;} \vspace{-.5mm}
\item [(e)] {\it $D\oplus D'$ where $D$ and $D'$ are division rings; or}
\vspace{-.5mm}
\item [(f)] {\it ${\Bbb Z}_3\oplus B\oplus D$ where $B$ is a Boolean ring and $D$ is a division ring.}
\vspace{-.5mm}
\item [(g)] {\it ${\Bbb Z}_3\oplus {\Bbb Z}_3\oplus D$;}\vspace{-.5mm}
\item [(h)] {\it ${\Bbb Z}_3\oplus {\Bbb Z}_3\oplus {\Bbb Z}_3\oplus {\Bbb Z}_3$;}
\vspace{-.5mm}
\item [(i)] {\it ${\Bbb Z}_3\oplus {\Bbb Z}_3\oplus {\Bbb Z}_3\oplus B$ where $B$ is adivision ring;}\vspace{-.5mm}
\item [(j)] {\it ${\Bbb Z}_3\oplus {\Bbb Z}_3\oplus {\Bbb Z}_3\oplus {\Bbb Z}_3$;}\end{enumerate}

In Case (b). If $0,\pm 1\neq x\in D$, then $(0,x)\not\in
U(R)\bigcup Id(R)\bigcup -Id(R)$. This forces $D\cong {\Bbb Z}_2,
{\Bbb Z}_3$. Hence, $(b)$ forces $R$ is in $(1)$ or $(3)$. $(c)$
does not occur. $(e)$ forces $D,D'\cong {\Bbb Z_2}$ or ${\Bbb
Z}_3$. Hence, $R$ is in $(1)-(3)$. $(f)$ does not occur except
$D\cong {\Bbb Z}_2$. Thus, $R$ is in $(1)-(3)$. Cases (g)-(j) do
not occur as $(1,-1,0), (1,-1,0,0)\not\in I(R)\bigcup
-Id(R)\bigcup N(R)$, as desired.

$\Longleftarrow$ This is clear.\hfill$\Box$

\vskip4mm \hspace{-1.8em} {\bf Theorem 4.3.}\ \ {\it Let $R$ be an
abelian ring. Then $R$ consists entirely of very idempotents,
units, and nilpotent elements if and only if $R$ is isomorphic to
one of the following:}\vspace{-.5mm}
\begin{enumerate}
\item [(1)] {\it ${\Bbb Z}_3$;}
\vspace{-.5mm}
\item [(2)] {\it a Boolean ring;}
\vspace{-.5mm}
\item [(3)] {\it ${\Bbb Z}_3\oplus {\Bbb Z}_3$;}
\vspace{-.5mm}
\item [(4)] {\it
${\Bbb Z}_3\oplus B$ where $B$ is a Boolean ring; or}
\item [(5)] {\it local ring with nil Jacobson radical.}
\end{enumerate} \vspace{-.5mm} {\it Proof.}\ \ $\Longrightarrow$
Case I. $R$ is indecomposable. Then $R=U(R)\bigcup N(R)$. This
shows that $R$ is local. Let $x\in J(R)$, then $x\in N(R)$, and so
$J(R)$ is nil.

Case II. $R$ is decomposable. In view of Lemma 4.2, $R$ is
isomorphic to one of the following: \vspace{-.5mm}
\begin{enumerate}
\item [(1)] {\it a Boolean ring;}
\vspace{-.5mm}
\item [(2)] {\it ${\Bbb Z}_3\times {\Bbb Z}_3$;}
\vspace{-.5mm}
\item [(3)] {\it
${\Bbb Z}_3\oplus B$ where $B$ is a Boolean ring.}
\end{enumerate} This shows that $R$ is isomorphic to one of
$(1)-(5)$, as desired.

$\Longleftarrow$ This is obvious.\hfill$\Box$

\vskip4mm \hspace{-1.8em} {\bf Lemma 4.4.}\ \ {\it Let $R$ be any
ring that consists entirely of very idempotents, units and
nilpotent elements. Then $eRe$ is a division ring for any
noncentral idempotent $e\in R$.}\vskip2mm\hspace{-1.8em} {\it
Proof.}\ \ Let $e\in R$ be a noncentral idempotent, and let
$f=1-e$. Then $R\cong \left(
\begin{array}{cc}
eRe&eRf\\
fRe&fRf
\end{array}
\right)$. The subring $\left(
\begin{array}{cc}
eRe&0\\
0&fRf
\end{array}
\right)$ consists entirely of very idempotents, units and
nilpotent elements. That is, $eRe\times fRf$ consists entirely of
very idempotents, units and nilpotent elements. Set $A=eRe$ and
$B=fRf$. Similarly to Lemma 4.2, $A=U(A)\bigcup Id(A)\bigcup
-Id(A)$. In view of Lemma 4.1, $A$ is isomorphic to one of the
following:
\begin{enumerate}
\item [(1)] {\it ${\Bbb Z}_3$;}
\vspace{-.5mm} \item [(2)] {\it a Boolean ring;} \vspace{-.5mm}
\item [(3)] {\it a division ring;}
\vspace{-.5mm}
\item [(4)] {\it ${\Bbb Z}_3\oplus B$ where $B$ is a Boolean ring.}
\end{enumerate}
That is, $A$ is a division ring or a ring in which every element
is very idempotent. Suppose that $eRe$ is not a division ring.
Then $eRe$ must contains a nontrivial idempotent, say $a\in R$.
Let $b=e-a$. Let $x\in eRf$ and $y\in fRe$. Choose
$$X_1=\left(
\begin{array}{cc}
a&x\\
0&0
\end{array}
\right),X_2=\left(
\begin{array}{cc}
b&x\\
0&0
\end{array}
\right),Y_1=\left(
\begin{array}{cc}
a&0\\
y&0
\end{array}
\right),Y_2=\left(
\begin{array}{cc}
b&0\\
y&0
\end{array}
\right).$$ Then $X_1,X_2,Y_1,Y_2$ are not invertible. As $a,b\in
eRe$ is nontrivial idempotents, we see that $X_1,X_2,Y_1,Y_2$ are
all not nilpotent matrices. This shows that $X_1$ and $X_2$ are
both very idempotents. It follows that $X_1=\pm X_2^2$ or
$X_2^2=\pm X_2$. As $x\in eRf, y\in fRe$, we have $ex=x$ and
$fy=y$.

Case I. $X_1=X_1^2, X_2=X_2^2$. Then $ax=x, bx=x$, and so
$x=ex=2x$; hence, $x=0$.

Case II. $X_1=X_1^2, X_2=-X_2^2$. Then $ax=x, bx=-x$, and so
$x=ex=0$.

Case III. $X_1=-X_1^2, X_2=X_2^2$. Then $ax=-x, bx=x$, and so
$x=ex=0$.

Case IV. $X_1=-X_1^2, X_2=-X_2^2$. Then $a=-a^2, ax=-x, b=-b^2$
and $bx=-x$. Hence, $(e-a)x=-x$, and so $x=ex=-2x$, hence, $3x=0$.
As $a\in R$ is an idempotent, we see that $a=a^2$, hence, $a=-a$,
and so $2a=0$. It follows that $x=-ax=(2a)x-(3x)a=0$.

Thus, $x=0$ in any case. We infer that $eRf=0$. Likewise, $fRe=0$.
Hence, $e\in R$ is central, an absurd. This completes the
proof.\hfill$\Box$

\vskip4mm \hspace{-1.8em} {\bf Lemma 4.5.}\ \ {\it Let $R$ be any
ring that consists entirely of very idempotents, units and
nilpotent elements. Then $eRe$ is isomorphic to ${\Bbb Z}/2{\Bbb
Z}$ or ${\Bbb Z}/3{\Bbb Z}$ for any noncentral idempotent $e\in
R$.}\vskip2mm\hspace{-1.8em} {\it Proof.}\ \ Let $e\in R$ be a
noncentral idempotent. In view of Lemma 4.4, $eRe$ is a division
ring. Set $f=1-e$. For any $u\in eRe$, we assume that $u\neq 0, e,
-e$, then the matrix
$$X= \left(
\begin{array}{cc}
u&0\\
0&0
\end{array}
\right)\in \left(
\begin{array}{cc}
eRe&eRf\\
fRe&fRf
\end{array}
\right)$$ is not be a unit, a very idempotent, or a nilpotent
element. This gives a contradiction. Therefore $u=0,e$ or $-e$, as
desired.\hfill$\Box$

\vskip4mm Recall that a ring $R$ is semiprime if it has no nonzero
nilpotent ideals. Furthermore, we derive

\vskip4mm \hspace{-1.8em} {\bf Theorem 4.6.}\ \ {\it Let $R$ be
any nonabelian ring that consists entirely of units, very
idempotents, and nilpotent elements. If $R$ is semiprime, then it
is isomorphic to $M_2({\Bbb Z}_2)$ or $M_2({\Bbb
Z}_3)$.}\vskip2mm\hspace{-1.8em} {\it Proof.}\ \ Suppose that $R$
is semiprime. In view of Lemma 4.4, $eRe$ is a division ring for
any noncentral idempotent $e\in R$. It follows by ~\cite[Lemma
21]{Du} that $R$ is isomorphic to $M_2(D)$ for a division ring
$D$. Choose $E_{11}=\left(
\begin{array}{cc}
1&0\\
0&0
\end{array}
\right)\in M_2(D)$. Then $E_{11}$ is a noncentral idempotent.
According to Lemma 4.5, $R\cong {\Bbb Z}/2{\Bbb Z}$ or ${\Bbb
Z}/3{\Bbb Z}$, as asserted.\hfill$\Box$

\vskip4mm Recall that a ring $R$ is a NJ-ring provided that for
any $a\in R$, either $a\in R$ is regular or $1-a\in R$ is a
unit~\cite{N}. Clearly, all rings in which every elements consists
entirely of units, very idempotents, and nilpotent elements are
NJ-rings.

\vskip4mm \hspace{-1.8em} {\bf Theorem 4.7.}\ \ {\it Let $R$ be
any nonabelian ring that consists entirely of very idempotents,
units and nilpotent elements. If $R$ is not semiprime, then it is
isomorphic to the ring of a Morita context with zero pairings
where the underlying rings are ${\Bbb Z}_2$ or ${\Bbb
Z}_3$.}\vskip2mm\hspace{-1.8em} {\it Proof.}\ \ Suppose that $R$
is not semiprime. Clearly, $R$ is a NJ-ring. In view of
~\cite[Theorem 2]{N}, $R$ must be a regular ring, a local ring or
isomorphic to the ring of a Morita context with zero pairings
where the underlying rings are both division ring. If $R$ is
regular, it is semiprime, a contradiction. If $R$ is local, it is
abelian, a contradiction. Therefore, $R$ is isomorphic to the ring
of a Morita context $T=(A,B,M,N,\varphi,\psi)$ with zero pairings
$\varphi,\psi$ where the underlying rings are division rings $A$
and $B$. Choose $E=\left(
\begin{array}{cc}
1_A&0\\
0&0
\end{array}
\right)\in T$. Then $E\in T$ is a noncentral idempotent. In light
of Lemma 4.5, $A\cong ETE\cong {\Bbb Z}_2$ or ${\Bbb Z}_3$.
Likewise, $B\cong  {\Bbb Z}_2$ or ${\Bbb Z}_3$. This completes the
proof.\hfill$\Box$

\vskip4mm With these information we completely determine the
structure of rings that consist entirely of very idempotents,
units and nilpotent elements.

\vskip4mm \hspace{-1.8em} {\bf Theorem 4.8.}\ \ {\it Let $R$ be a
ring. Then $R$ consists entirely of very idempotents, units and
nilpotent elements if and only if $R$ is isomorphic to one of the
following:}\vspace{-.5mm}
\begin{enumerate}
\item [(1)] {\it a Boolean ring;}
\vspace{-.5mm}
\item [(2)] {\it ${\Bbb Z}_3\times {\Bbb Z}_3$;}
\vspace{-.5mm}
\item [(3)] {\it
${\Bbb Z}_3\oplus B$ where $B$ is a Boolean ring;}
\item [(4)] {\it local ring with a nil Jacobson radical;}
\vspace{-.5mm}
\item [(5)] {\it $M_2\big({\Bbb Z}_2\big)$ or $M_2\big({\Bbb Z}_3\big)$; or} \vspace{-.5mm}
\item [(6)] {\it the ring of a Morita context
with zero pairings where the underlying rings are ${\Bbb Z}_2$ or
${\Bbb Z}_3$.}
\end{enumerate} \vspace{-.5mm} {\it Proof.}\ \ $\Longrightarrow$
This is obvious by Theorem 4.3, Theorem 4.6 and Theorem 4.7.

$\Longleftarrow$ Cases (1)-(4) are easy. Case (5)-(6) are verified
by checking all possible (generalized) matrices over ${\Bbb Z}_2$
and ${\Bbb Z}_3$.\hfill$\Box$

\end{document}